\documentclass{article}

\usepackage{amsmath,amssymb,enumerate,bbm,calc,capt-of,ifthen,mathrsfs}
\usepackage{epsfig}
\usepackage{url}

\newtheorem{cntr}{ERROR! SHOULD NOT USE THIS}

\newcommand{\tmop}[1]{\operatorname{#1}}

\newtheorem{varremark}[cntr]{Remark}

\newcommand{\nin}{\not\in}

\newtheorem{varnote}{Note}

\newenvironment{proof}{
  \noindent\textbf{Proof.}\ }{\hspace*{\fill}
  \begin{math}\Box\end{math}\medskip}

\newenvironment{proof*}[1]{
  \noindent\textbf{#1\ }}{\hspace*{\fill}
  \begin{math}\Box\end{math}\medskip}

\newtheorem{theorem}[cntr]{Theorem}

\newtheorem{algo}{Algorithm}

\newcommand{\sdev}[0]{\sigma}
\newcommand{\ndom}[0]{N}
\newcommand{\nran}[0]{n}

\newcommand{\RtwoN}[0]{\mathbb{R}^{\ndom} \times \mathbb{R}^{\ndom}}

\newcommand{\Rn}[0]{{\mathbb{R}^{\ndom}}}
\newcommand{\Pout}[2]{{\mathcal{O}_{#1}^{#2} }}

\newcommand{\kmax}[0]{k_{\tmop{max}}}
\newcommand{\vmax}[0]{v_{\tmop{max}}}

\newcommand{\Tstep}[0]{T_{\tmop{step}}}

\newcommand{\Lb}{{L}}
\newcommand{\wb}{{w}}
\newcommand{\kb}{{k_{b}}}

\newcommand{\vx}[0]{\vec{x}}

\newcommand{\vk}[0]{\vec{k}}

\newcommand{\vu}{{\vec{u}}}
\newcommand{\vua}{{\vec{u}_{d}}}
\newcommand{\vd}{{\vec{d}}}
\newcommand{\HH}{{\mathcal{H}}}
\newcommand{\HHb}{{\mathcal{H}_{b}}}
\newcommand{\U}[1]{{e^{\HH #1}}}
\newcommand{\Ub}[1]{{e^{\HHb #1}}}

\newcommand{\norm}[2]{\left\| #1 \right \|_{#2}}
\newcommand{\abs}[1]{\left| #1 \right|}
\newcommand{\absSmall}[1]{| #1 |}

\newcommand{\IP}[2]{\left\langle #1 | #2 \right\rangle}

\newcommand{\RR}{{\mathbb{R}}}
\newcommand{\CC}{{\mathbb{C}}}

\newcommand{\turnaround}{{Z}}
\newcommand{\aturnaround}{{Z}}

\numberwithin{algo}{section}
\numberwithin{equation}{section}

\newcommand{\comment}[1]{}

\begin{document}

\title{A stable absorbing boundary layer for anisotropic waves}

\author{A. Soffer and C. Stucchio}

\maketitle

\begin{abstract}
  For some anisotropic wave models, the PML (perfectly matched layer) method of open boundaries can become polynomially or exponentially unstable in time. In this work we present a new method of open boundaries, the phase space filter, which is stable for all wave equations.

  Outgoing waves can be characterized as waves located near the boundary of the computational domain with group velocities pointing outward. The phase space filtering algorithm consists of applying a filter to the solution that removes outgoing waves only.

 The method presented here is a simplified version of the original phase space filter, originally described in \cite{us:TDPSFjcp} for the Schr\"odinger equation. We apply this method to anisotropic wave models for which the PML is unstable, namely the Euler equations (linearized about a constant jet flow) and Maxwell's equations in an anisotropic medium. Stability of the phase space filter is proved.
\end{abstract}

\section{Introduction}

We consider the solution of linear wave equations, of the form:

\begin{equation}
  \label{eq:waveEq}
  \vu_{t}(\vx,t) = \HH \vu(\vx,t)
\end{equation}
Here, $\vu: \RR^{\ndom} \rightarrow \RR^{\nran}$ (or $\RR^{\ndom} \rightarrow \CC^{\nran}$) and $\HH$ is a skew-adjoint linear differential operator. This is a prototypical example of a linear wave equation. Such equations admit linear waves of the form $e^{i (\vk \cdot \vx-\omega_{j}(\vk)t)}\vd_{j}(k)$ with $\vd_{j}(\vk)$ the $j$'th eigenvector of $\HH$ and $\omega_{j}(\vk)$ the $j$-th eigenvalue of $\HH$ in the frequency domain.

A great many systems of practical importance can be recast in such a form, including the wave equation, Schr\"odingers equation, the linearized Euler equations, Maxwell's equations and others. Very often, the equations are too complicated to solve exactly, so numerical approximations must be used. If we are only interested in the solution of $\vu(x,t)$ in some finite region $B \subset \Rn$, it simplifies the computations to solve \eqref{eq:waveEq} only on this region. Of course, steps must be taken to prevent spurious reflection from the artificial boundary.

Exact non-reflecting boundary conditions (NRBCs) can be constructed for many wave equations \cite{MR1913093,MR596431,MR0471386,MR0436612,MR517938,MR2032866,jian:thesis}. Such methods are accurate, but non-local in time and space, and algorithms are sensitive to the shape of the boundary and the particular interior solver. They also prevent the use of Fourier spectral methods for solving the interior problem, since these methods have periodic or Dirichlet boundaries built in.

The Perfectly Matched Layer (PML) is another approach to the problem \cite{MR1294924,MR1412240}. The PML consists of the complex change of coordinates, $\vx \mapsto \vx+i\sigma(\vx)$ (for a suitably chosen $\sigma(\vx)$) where $\sigma(\vx)=0$ for $\vx \in B$ (this technique was used in scattering theory for a long time \cite{yajima:ACStarkExteriorScaling,simon:resonanceReview}). Inside $B$ the dynamics of \eqref{eq:waveEq} are unchanged, while on $B^{C}$ plane waves are turned into exponentially decaying modes. The domain of computation is taken to be a region $B_{c} \supset B$ which is $B$ surrounded by an absorbing layer. The layer dissipates outgoing waves and introduces minimal spurious reflections. Unlike NRBCs, the PML is local in space and time and compatible with fast spectral methods. For this reason, the PML is the method of choice in modern wave propagation.

However, the PML has a fatal flaw which limits its use in some cases. For certain types of anisotropic wave, the PML can become exponentially unstable in time regardless of the particular method used (the instability exists on the level of the PDE). This was first noticed in \cite{hu:unstablePML} for the Euler equations linearized about a jet flow (see also \cite{MR1294924,becache:unstablePML,pmlShallowWater}). In \cite{becache:unstablePML}, a simple geometrical criterion\footnote{There are some other conditions necessary for the proof to apply, but it is convincingly argued that this is the fundamental criteria.} was provided for a PML to be unstable. Consider a boundary at $\vx_{1}=L$ and a dispersion relation $\omega(\vk)$: if  $\vk_{1} > 0$ while the group velocity $\partial_{k_{1}} \omega(\vk)<0$ for some $\vk$, then the PML is exponentially unstable.

In \cite{us:TDPSFrigorous, us:TDPSFjcp} a new approach to the problem of open boundaries was introduced, the Time Dependent Phase Space Filter (TDPSF). Phase space analysis is a major tool in modern scattering theory \cite{MR818063,MR898052,MR1233895}. The key idea in this approach is that certain regions of quantum phase space (the set of points $\{ (\vx,\vk) \in \RtwoN \}$, where $\vx$ represents a position and $\vk$ a spatial frequency) consist solely of outgoing waves, whereas other regions have more complicated interactions. The philosophy of the TDPSF is to identify parts of the solution $\vu(\vx,t)$ localized in the outgoing regions of phase space, and filter them from the solution before they reach the boundary. After these waves are filtered, $\vu(\vx,t)$ is not approaching the boundary, and therefore boundary conditions don't matter.

In \cite{us:multiscale,us:TDPSFrigorous,us:TDPSFjcp}, this approach is used to construct open boundaries for the Schr\"odinger equation. The phase space projections there are based on the Gaussian windowed Fourier transform \cite{MR966733}. In this paper, we extend the work of \cite{us:multiscale,us:TDPSFrigorous,us:TDPSFjcp} to more general wave equations including the linearized Euler equation. We also simplify the method significantly, replacing the windowed Fourier transform by standard phase space projections of the form $\chi(\vx) P(k) \chi(\vx)$ (this approach was also taken in \cite{us:multiscale}). We begin by briefly reviewing the dynamics of linear waves (Section \ref{sec:dynamicsOfWaves}). In Section \ref{sec:AlgorithmSection} we present the phase space filtering algorithm for general systems of the form \eqref{eq:waveEq}. In Section \ref{sec:Examples}, we show some numerical examples. In Section \ref{sec:lowFrequency}, we briefly discuss a method of filtering outgoing waves with frequencies too small to resolve on an absorbing boundary layer of reasonable width.

\subsection{Dynamics of Waves}
\label{sec:dynamicsOfWaves}

For concreteness and to pin down our notation, we review a few facts about wave propagation. Recall that we defined $\omega_{j}(\vk)$ as the $k$'th branch of the dispersion relation, and $\vd_{j}(\vk)$ as the corresponding normalized eigenvector. With this notation, $e^{i (\vk \cdot \vx-\omega_{j}(k) t)} d_{j}(\vk)$ is a standard plane wave solution of \eqref{eq:waveEq}.

Let us define the matrix $D$ to be the (unitary) matrix with $j$'th row $\vd_{j}(k)$. This matrix can be used to diagonalize $\HH$, i.e:
\begin{equation}
  \label{eq:diagonalizeH}
  \HH = D^{\dagger} \left[
    \begin{array}{lll}
      i\omega_{1}(\vk) & \ldots & 0\\
      \ldots & i\omega_{j}(\vk) &  \ldots\\
      0 & \ldots & i\omega_{\nran}(\vk)\\
    \end{array}
    \right] D
\end{equation}

Let $\U{t}$ denote the propagation operator for \eqref{eq:waveEq}, i.e. the operator mapping $\vu(x,0)$ to $\vu(x,t)$. In the frequency domain, the propagator can be written as:
\begin{equation}
  \label{eq:propagatorDef}
  \U{t} = D^{\dagger} \exp \left( \left[
    \begin{array}{lll}
      i\omega_{1}(\vk) & \ldots & 0\\
      \ldots & i\omega_{j}(\vk) &  \ldots\\
      0 & \ldots & i\omega_{\nran}(\vk)\\
    \end{array}
    \right] t \right) D
\end{equation}

Now consider an initial condition $\vu_{0}(x)$ localized in frequency about the point $\vk_{0}$ and in position about the point $\vx_{0}$. Consider the solution $\vu(x,t)$ with $\vu(x,t=0)=\vu_{0}(x)$. This solution can be written in the frequency domain as:
\begin{multline}
  \label{eq:propagatorFrequencyDomain}
  \widehat{\vu}(\vk,t) = \U{t} \vu_{0}(x) \widehat{\vu}_{0}(\vk)  \\
  = D^{\dagger} \exp \left( \left[
    \begin{array}{lll}
      i\omega_{1}(\vk) & \ldots & 0\\
      \ldots & i\omega_{j}(\vk) &  \ldots\\
      0 & \ldots & i\omega_{\nran}(\vk)\\
    \end{array}
    \right] t \right) \left(
    \begin{array}{l}
      u^{1}(\vk)\\
      \ldots\\
      u^{\nran}(\vk)
    \end{array}
  \right)
\end{multline}
where $u^{j}(\vk)$ is the projection of $\widehat{\vu}_{0}(\vk)$ onto $\vd_{j}(\vk)$ and $\widehat{f}(\vk)$ denotes the Fourier transform of $f(\vx)$.

The $j$-th component corresponds to a superposition of plane waves propagating with velocity $\nabla_{k}\omega_{j}(\vk)$. If we consider only regions of $\vk$-space in which $\partial_{k_{1}} \omega_{j}(\vk) > 0$, then we are considering only waves with a rightward moving component. The same can be said about other directions. It is this property combined with phase space localization techniques which we will use to filter outgoing waves.

\section{Method}
\label{sec:AlgorithmSection}

The methodology of the phase space filter is rather simple. First, given the operator $\HH$, we find the generalized eigenfunctions and the dispersion relation. Suppose that the dispersion relations $\omega_{j}(\vk)$ and $\vd_{j}(\vk)$ are given. If this is the case, then $\vd_{j}(\vk)$ is a plane wave propagating with group velocity $\nabla_{k} \omega^{j}(\vk)$.

\subsection{The Propagation Algorithm}

\label{sec:propagationAlgorithm}

The propagation algorithm is as follows. We solve \eqref{eq:waveEq} on the region $[-\Lb-\wb,\Lb+\wb]^{\ndom}$ with any accurate interior solver. At this point, we assume the initial condition $\vu_{0}(x)$ is approximately band-limited. We assume that the mass of $\widehat{\vu}_{0}(\vk)$ for $\abs{\vk} \geq \kmax$ is small. This assumption is common: these are simply frequencies too high to resolve with the computational grid.

We also assume that the frequency content of $\widehat{\vu}_{0}(\vk)$ is not localized near the regions where the group velocity turns around. I.e., let $\turnaround=\{ \vk : \exists i,j \partial_{k_{i}} \omega_{j}(\vk) =0\}$ be the set of points where the group velocity vanishes in some direction. Let $\aturnaround$ be all points a distance $\kb$ from $Z$. We \emph{assume} that for some small $\delta$,
\begin{equation}
  \int_{\aturnaround} \abs{\widehat{\vu}_{0}(\vk)}^{2} d\vk \leq \delta
\end{equation}
This requirement is necessary due to the Heisenberg uncertainty principle. Very roughly, the Heisenberg uncertainty principle says that if the width of the filter region is $\wb$, then we cannot localize on any region of $\vk$-space smaller than $\wb^{-1}$ (up to constants). To localize only outgoing waves, we would need to isolate waves with group velocity pointing outward, while not localizing on waves with group velocity pointing inward. Localizing this accurately on waves inside $\aturnaround$ would be impossible, due to the uncertainty principle. This difficulty can be resolved at logarithmic cost, and this is discussed briefly in Section \ref{sec:lowFrequency}.

The propagation algorithm is simple. Fix a time $\Tstep \leq \wb / 3 \vmax$, with
\begin{equation}
\vmax= \sup_{j} \sup_{\abs{\vk} \leq \kmax} \abs{\nabla_{k} \omega_{j}(\vk)}
\end{equation}
being the the largest group velocity relevant to the problem. This criterion ensures that waves cannot cross the buffer region in a time interval shorter than $\Tstep$.

On the time intervals $[0,\Tstep]$, $[\Tstep, 2\Tstep]$, $\ldots$, we solve \eqref{eq:waveEq} with the interior propagator. At times $\Tstep$, $2 \Tstep$, $\ldots$, we apply outgoing wave filters $\Pout{j}{\pm}$ (defined in the next section) to the regions $[-\Lb-\wb,\Lb+\wb]^{\ndom} \setminus [-\Lb,\Lb]^{\ndom}$. After the application of this filter, all waves which would have reached the boundary before a time $\Tstep$ have been removed. Thus, no waves actually reach the boundary, and the boundary conditions are irrelevant.

\begin{algo}
  \label{algo:propagation}
  TDPSF Propagation Algorithm

  Given an initial condition $\vu_{0}(x)$, this algorithm calculates $\vua(x,t)$.

  { \bf Input: }
  \begin{itemize}
  \item The dispersion relations and diagonalizing matrices, $\omega_{j}(\vk)$ and $D$.
  \item $\Ub{t}$, a propagator that accurately solves the interior problem.
  \item $\kmax$, the maximal frequency of the problem.
  \end{itemize}

  \begin{enumerate}
  \item Define $\Tstep=\wb/3v_{max}$.
  \item Define the approximate solution $\vua(x,t)$ recursively.
    \begin{subequations}
      At times which are an integer multiple of $\Tstep$, we filter off the outgoing waves:
      \begin{equation}
        \vua(x,(m+1)\Tstep) = \left[ \prod_{j=1}^{\ndom} (1-\Pout{j}{+})(1-\Pout{k}{-}) \right] \Ub{\Tstep} \vua(x,m \Tstep)
      \end{equation}
      The outgoing wave filters are computed using Algorithm \ref{algo:outgoingWaveFilter}, described in Section \ref{sec:boundaryFilterDetails}.

      For other times, we use the given interior propagator:
      \begin{equation}
        \vua(x,m \Tstep+\tau) = \Ub{\tau} \vua(x,m \Tstep) \textrm{ for } \tau \in [0,\Tstep)
      \end{equation}
      \begin{equation}
        \vua(x,0) = \vu_{0}(x)
      \end{equation}
    \end{subequations}
  \end{enumerate}
\end{algo}

It now remains to construct the outgoing wave filters, $\Pout{j}{\pm}$.

\subsection{Construction of the boundary filter}
\label{sec:boundaryFilterDetails}

Consider a fixed boundary region, say the boundary at $\vx_{1}=L$. For a frequency $\vk$, if $\partial_{1} \omega_{j}(\vk) > 0$, then waves with frequency $\vk$ are outgoing at this boundary. The outgoing region of phase space at the right boundary is therefore
\begin{equation}
  \{ (\vx,\vk) \in \RtwoN : \vx_{1} > L \textrm{ and } \partial_{1} \omega_{j}(\vk) > 0\}.
\end{equation}

We now construct a projection onto this region. The Heisenberg uncertainty principle makes an exact projection impossible, but we will do the best it allows.

Extend the box a width $\wb$, to be specified shortly. Define the function
\begin{equation}
  \label{eq:defOfChi}
  \chi^{\pm}_{j}(x) = \left(\frac{2}{\sigma \sqrt{\pi}}\right)^{\ndom} e^{-\vx^{2}/\sdev^{2}} \star I_{j}(x)
\end{equation}
where $I_{j}(x)=1$ for $x_{j} \in [\pm (\Lb+\wb/3),\pm (\Lb+2\wb/3)]$ and $x_{k} \in  [-\Lb-2\wb/3,\Lb+2\wb/3]$ (for $k \neq j$), and $0$ elsewhere. The parameter $\sigma=O(\wb/ \ln(\delta^{-1})^{1/2})$; a precise bound is given in \eqref{eq:upperBoundOnSigma}. This ensures that $\chi^{\pm}_{j}(x)<\delta$ for $x_{j} \nin [\pm \Lb, \pm (\Lb+\wb)$ or $x_{k} \nin [-\Lb-\wb,\Lb+\wb]$. This function is smooth, and well localized inside the buffer region on the $j$'th sides of the box.

The set $R_{j,l} = \{ \vk \in \Rn : \partial_{k_{j}} \omega_{l}(\vk) > 0 \}$ is the set of frequencies with the $k$'th branch of the group velocity pointing right. Due to the Heisenberg uncertainty principle, we cannot project onto this set precisely. However, we can approximately project onto most of these wave-vectors. Define
\begin{subequations}
  \begin{equation}
    R_{j,l,\delta} = \{ \vk \in R_{k} : d(\vk, R_{j,l}^{C}) > \kb \}
  \end{equation}
  The set of vectors within a distance $\kb$ of group velocity $\partial_{k_{j}} \omega_{l}(x)=0$ is a the set of group velocities with motion normal to side $j$ approximately equal to zero. The width $\kb$ is a buffer to ensure that the frequency spreading caused by our spatial localization operators does not cause an error larger than $\delta$. Given $\kb$, we must also choose $\sigma \geq O(\kb^{-1} \ln(\delta^{-1})^{1/2})$ to ensure that the spatial localization operators do not spread the frequency content of the solution past the buffer region of width $\kb$.

In particular, we must choose the buffers widths $\kb,\wb$ and standard deviation $\sigma$ to satisfy
\begin{multline}
  \label{eq:estimateOnkbsdev}
  \kb^{-1} \left(\ln(\delta^{-1}) + \ln\left( \frac{\wb^{2} \Lb^{\ndom-1} 2^{3\ndom} \sigma^{3\ndom} }{\pi^{3\ndom/2}} \right) \right)^{1/2} \\
  \leq \sigma
  \leq \frac{\wb}{
    \sqrt{\ln(\delta^{-1}) + \ndom \ln(2 \sigma \pi^{-1/2})  }
  }
\end{multline}
which ensures that spreading in frequency does not turn waves around (thus minimizing reflection). This is estimated in Section \ref{sec:choosingBufferWidths}. In particular, note that \eqref{eq:estimateOnkbsdev} implies the buffer width $\wb$ must be at least  $O( \kb^{-1} \ln(\delta^{-1}))$ (as expected from the Heisenberg uncertainty principle).

  The frequency projection operator is defined as:
  \begin{equation}
    P_{j,l,\delta}(k) = \left(\frac{2\sigma}{\sqrt{\pi}}\right)^{\ndom} e^{-k^{2} \sigma^{2}} \star
    \left[
    \begin{array}{lll}
      1_{R_{j,l,\delta}} (\vk)& \ldots & 0\\
      \ldots & 1_{R_{j,l,\delta}}(\vk) &  \ldots\\
      0 & \ldots & 1_{R_{j,n,\delta}}(\vk)\\
    \end{array}
    \right]
  \end{equation}
\end{subequations}
Thus, the operator $P_{j,l,\delta}(k)$ is a smooth projection (in the basis of eigenvectors of $\HH$) onto wave-vectors propagating rightward. Conjugation of $P(\vk)$ by $D$ moves the projection into the domain of frequency vectors. Finally, we define the operator:
\begin{equation}
  \label{eq:defOfPout}
  \Pout{1}{+} = \chi_{1}^{+}(x) D^{\dagger}  P(\vk) D \chi_{1}^{+}(x)
\end{equation}

This operator both localizes in the buffer region at $x=L$, and projects onto waves with group velocity pointing to the right. This operator can be computed efficiently and with spectral accuracy as follows:

\begin{algo}
  \label{algo:outgoingWaveFilter}
  Outgoing Wave Filter Algorithm

  This algorithm applies the outgoing wave filters, i.e. it numerically approximates $(1-\Pout{j}{\pm})\vu(x)$.

  { \bf Input: }
  \begin{itemize}
  \item The dispersion relations and diagonalizing matrices, $\omega_{j}(\vk)$ and $D$.
  \item A function $\vu(x)$.
  \item Indices $(j, \cdot)$ with $j \in 1 \ldots \nran$ and $\cdot \in \{+,-\}$.
  \end{itemize}
  We assume that the operator $D^{\dagger} P(\vk) D$ is precomputed.
  \begin{enumerate}
  \item Compute the function $\chi_{j}^{\pm}(x) \vu(x)$, and take it's Fast Fourier Transform.
  \item Apply $D^{\dagger} P(\vk) D$ to the Fast Fourier transform of $\chi_{j}^{\pm}(x) \vu(x)$.
  \item Apply the inverse Fast Fourier Transform, and multiply the result by $\chi_{j}^{\pm}(x)$.
  \item Subtract the result from $\vu(x)$, and return the result.
  \end{enumerate}
\end{algo}

\subsubsection{Choosing the Buffer Widths}
\label{sec:choosingBufferWidths}

The buffer widths, $\wb$ and $\kb$ must be chosen carefully for this algorithm to work. There are two competing concerns, namely width of the buffer region and frequency spreading, which we must address at this point.

If $\wb$ is too small, then the spatial localization functions $\chi_{j}^{\pm}(x)$ will become rougher. But roughness in $\chi_{j}^{\pm}(x)$ will spread the frequencies of the solution around. In the frequency domain, multiplication by $\chi_{j}^{\pm}(x)$ corresponds to convolution, and behaves much like a diffusion operator in the $k$ variable. The danger is that the diffusion in $k$ might move mass from the region with positive group velocity to the region with negative group velocity.

However, a larger $\wb$ will increase the size of the computational box, and therefore the computational complexity of the method. Thus, it is desirable to take $\wb$ as small as possible.

We analyze this as follows. Multiplication by $\chi_{j}(x)$ corresponds (in the $k$-domain to convolution with\footnote{We have suppressed the $\textrm{sinc}(...)$-factor corresponding to the Fourier transform of $I_{j}(x)$, just for simplicity. This is an overestimate.} $(\wb \Lb^{\ndom-1}/3)(2^{\ndom}\sigma^{\ndom}\pi^{-\ndom/2})e^{-k^{2} \sigma^{2}}$. Thus, we can approximate the frequency domain operations by:
\begin{multline*}
  (\wb \Lb^{\ndom-1}/3)(2^{\ndom}\sigma^{\ndom}\pi^{-\ndom/2})e^{-k^{2} \sigma^{2}} \star (\wb \Lb^{\ndom-1}/3)(2^{\ndom}\sigma^{\ndom}\pi^{-\ndom/2})e^{-k^{2} \sigma^{2}} \star\\
  = \wb^{2} \Lb^{\ndom-1} 2^{2\ndom} \sigma^{2\ndom} \pi^{-\ndom} e^{-k^{2} \sigma^{2}/2} \star
\end{multline*}
The frequency domain operators $P_{j,l,\delta}(\vk)$ are localized on the regions $R_{j,k,\delta}$, also with a Gaussian tail of the form $(2\sigma \pi^{-1/2})^{\ndom} e^{-k^{2} \sigma^{2}}$.

The region $R_{j,l,\delta}$ is separated from the region containing incoming waves by a buffer of width $\kb$. We want to make sure that the spreading in frequency is small past this buffer of width $\kb$, i.e.:
\begin{equation*}
  \wb^{2} \Lb^{\ndom-1} 2^{2\ndom} \sigma^{2\ndom} \pi^{-\ndom} e^{-k^{2} \sigma^{2}/2} \star (2\sigma \pi^{-1/2})^{\ndom} e^{-k^{2} \sigma^{2}} \Bigg|_{\abs{\vk}=\kb} \leq \delta
\end{equation*}
This can be guaranteed by:
\begin{equation}
  \label{eq:LowerBoundOnSigma}
  \kb^{-1} \left(\ln(\delta^{-1}) + \ln\left( \frac{\wb^{2} \Lb^{\ndom-1} 2^{3\ndom} \sigma^{3\ndom} }{\pi^{3\ndom/2}} \right) \right)^{1/2} \leq \sigma
\end{equation}
On the other hand, if we take $\sigma$ too large, then the tails of $\chi_{j}^{\pm}(\vx)$ will enter the computational domain. To ensure that this is minimized, we want to make certain that $\chi_{j}^{\pm}(\vx) \leq \delta$ for $\vx \in [-\Lb,\Lb]^{\ndom}$. By examining the form of \eqref{eq:defOfChi}, we find that this can be accomplished if:
\begin{equation*}
  \left(\frac{2}{\sigma \sqrt{\pi}}\right)^{\ndom} e^{-\wb^{2}/\sigma^{2}} \leq \delta
\end{equation*}
or equivalently
\begin{equation}
  \label{eq:upperBoundOnSigma}
  \sigma \leq \frac{\wb}{
    \sqrt{\ln(\delta^{-1}) + \ndom \ln(2 \sigma \pi^{-1/2})  }
  }
\end{equation}
Additionally, to satisfy simultaneously \eqref{eq:LowerBoundOnSigma} and \eqref{eq:upperBoundOnSigma}, we must have that:
\begin{multline}
  \label{eq:constraintOnKbW}
\kb^{-1} \left(\ln(\delta^{-1}) + \ln\left( \frac{\wb^{2} \Lb^{\ndom-1} 2^{3\ndom} \sigma^{3\ndom} }{\pi^{3\ndom/2}} \right) \right)^{1/2} \\
  \leq \frac{\wb}{
    \sqrt{\ln(\delta^{-1}) + \ndom \ln(2 \sigma \pi^{-1/2})  }
  }
\end{multline}
This condition demands essentially that $\wb \geq C \kb^{-1}$, where $C=O(\ln(\delta^{-1}))$. This is precisely what we should expect based on the Heisenberg uncertainty principle.

\subsubsection{Powers of $2$}

As a practical matter, there is an additional constraint on the buffer width. The projections onto outward moving group velocities  (the $D^{\dagger}  P(\vk) D$ part of $\Pout{j}{\pm}$) are performed using an FFT. The FFT algorithm is fastest when performed on a grid of size $2^{m}$ in each dimension. For this reason, it is efficient to take $w=2^{m} \delta x$ with $m=\lceil \log_{2}(w_{min}/\delta x) \rceil$. Here, $w_{min}$ is minimal $w$ satisfying \eqref{eq:constraintOnKbW} and $\delta x$ is the lattice spacing in position.

\subsection{Stability}

\label{sec:stability}

The stability of the method is readily proved. The main reason that the algorithm is stable is simply that it is dissipative: the filtering operator $(1-\Pout{j}{\pm})$ has norm bounded by $1$ and can therefore not increase the norm of $\vu(x,t)$. Additionally, all dissipation occurs only at discrete instants of time, which minimizes interactions with the propagator. Using these ideas, it is straightforward to prove that the algorithm is stable, under the sole assumption that the interior propagator is stable.

\begin{theorem}
  \label{thm:stability}
  The Time Dependent Phase Space Filtering propagation algorithm is stable if the interior solver is. In particular, we have the estimate:
  \begin{equation}
    \label{eq:stability}
    \norm{\vua(x,t)}{L^{2}} \leq \norm{\vu_{0}(x)}{L^{2}}
  \end{equation}
\end{theorem}

\begin{proof}
  The idea of the proof is to show that the numerical solution operator at time $t$ can be written as the product of operators of norm $1$. Thus energy ($L^{2}$ norm) at time $t$ must be bounded by the initial energy.

  Define the operator $U$ by:
  \begin{equation*}
    U = \left[ \prod_{j=1}^{\ndom} (1-\Pout{j}{+})(1-\Pout{k}{-}) \right] \Ub{\Tstep},
  \end{equation*}
  At time $t=m \Tstep+\tau$ (with $\tau \in [0,\Tstep]$), we can write the numerical solution $\vua(x,t)$ as:
  \begin{equation*}
    \vua(x,t) = \Ub{\tau} U^{m} \vu_{0}(x)
  \end{equation*}
  By the self-adjointness of $\HHb$, $\norm{\Ub{t}}{} = 1$. If we can show that $\norm{1-\Pout{k}{\pm}}{} \leq 1$, then $\norm{U}{}=1$, implying:
  \begin{equation*}
    \norm{\vu_{b}(x,t)}{L^{2}} \leq \norm{\Ub{\tau}}{} \norm{U}{}^{m} \norm{\vu_{0}(x)}{L^{2}}
    \leq 1 \cdot 1^{m} \cdot \norm{\vu_{0}(x)}{L^{2}} = \norm{\vu_{0}(x)}{L^{2}}
  \end{equation*}
  and stability is proved.

  We first show that $\sigma(\Pout{j}{\pm}) \subseteq [0,1]$. Recall that $\Pout{j}{\pm}$ is defined as $\Pout{j}{\pm}=\chi_{j}^{\pm}(x) D^{\dagger} P(\vk) D \chi_{j}^{\pm}(x)$. Note that $D^{\dagger} P(\vk) D$ is diagonalizable (to $P(\vk)$) and each diagonal entry is contained in $[0,1]$ for each $\vk$. Therefore, $\norm{D^{\dagger} P(\vk) D}{} \leq 1$ and $\sigma(D^{\dagger} P(\vk) D) =\sigma(P(\vk)) \in [0,1]$. This implies that $\norm{\Pout{j}{\pm}}{} \leq 1$. Now write:
  \begin{multline*}
    \IP{ f}{\Pout{j}{\pm} f} =
    \IP{
      \chi_{j}^{\pm}(x) f(x)
    }{
      D^{\dagger} P(\vk) D
      \chi_{j}^{\pm}(x) f(x)
    } \\
    = \IP{
      D \chi_{j}^{\pm}(x) f(x)
    }{
      P(\vk) D
      \chi_{j}^{\pm}(x) f(x)
    } = \IP{g}{P(\vk) g}
  \end{multline*}
  where $g=D \chi_{j}^{\pm}(x) f(x)$. Since $P(\vk)$ is a positive matrix for each $\vk$, positivity follows. Since $\Pout{j}{\pm}$ is a positive operator with norm bounded by $1$, $\Pout{j}{\pm}$ has spectrum in $[0,1]$. The spectral mapping theorem implies that $\sigma(1-\Pout{j}{\pm}) \subseteq 1-[0,1] = [0,1]$, implying that $\norm{1-\Pout{j}{\pm}}{} \leq 1$, and stability is proved.
\end{proof}

\subsection{Tangential Waves}

The filter described in this section is not optimal. In particular, some outgoing waves that should be filtered are not, primarily waves which are outgoing, but nearly tangentially to the boundary.

While it is impossible to completely resolve this issue (due to the uncertainty principle), one can improve the situation by using a better phase space localization scheme. One can build phase space projections using framelets (localized wavepackets) fine tuned to the problem of interest. For the Schr\"odinger equation, canonical coherent states (the Gaussian Windowed Fourier transform) frame \cite{us:TDPSFjcp} is a natural choice. Other equations require different frames (e.g. curvelets or wave atoms for hyperbolic waves \cite{MR2165380}).

The cost of this approach is extra programmer time and (usually) an increase in computational complexity by a constant factor.

\section{Examples}
\label{sec:Examples}

\subsection{A warm-up: 1-dimensional Schr\"odinger equation}

Phase space filters were originally developed for the Schr\"odinger equation ($u(x,t)$ is a scalar field, and $\HH=i \Delta$). For the Schr\"odinger equation, the phase space filter takes a particularly simple form. There is no diagonalizing operator, and the group velocity is simply $k$. Thus, outgoing waves at the right boundary are simply waves of positive frequency (negative frequency on the left boundary).

In this case, $P(k)$ takes the simple form:
\begin{equation*}
  P(k) = \left( \frac{2\sigma}{\sqrt{\pi}}\right) e^{-k^{2}/\sigma^{2}} \star 1_{k > k_{min}} (k)
\end{equation*}

Thus, the filtering operator is simply
\begin{equation*}
  \Pout{1}{+} = \chi_{1}(x) \left( \frac{2\sigma}{\sqrt{\pi}}\right) \left[ e^{-k^{2}/\sigma^{2}} \star 1_{k > \kb} (k) \right]\chi_{1}(x)
\end{equation*}
(and similarly on the left side). For simplicity, we take $\kb=0$.

\begin{figure}
\setlength{\unitlength}{0.240900pt}
\ifx\plotpoint\undefined\newsavebox{\plotpoint}\fi
\sbox{\plotpoint}{\rule[-0.200pt]{0.400pt}{0.400pt}}%
\includegraphics[scale=0.5]{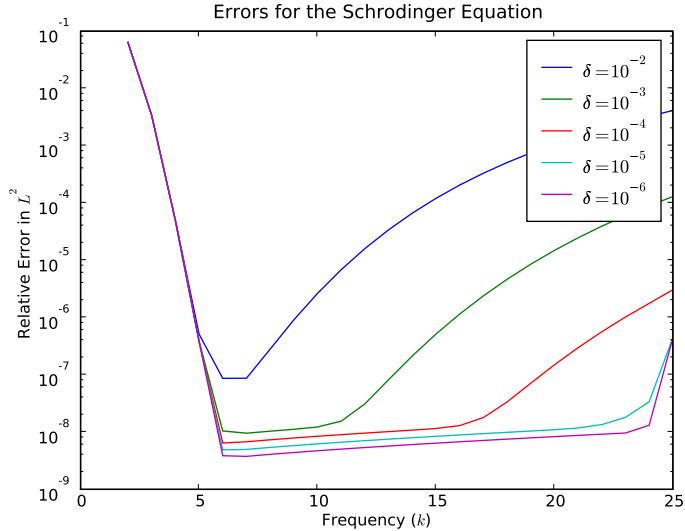}

\caption{The relative error for the 1-dimensional test of the Schrodinger equation, as a function of $k$ and the parameter $\delta$. }
\label{fig:schroErrors}
\end{figure}

We solved this example on a lattice of $1024$ pts, with $\delta x = 0.1$ and initial condition
\begin{equation*}
  u(x,t=0) = \frac{e^{i k x}}{2\sqrt{7}} e^{-x^{2}/2 \cdot 7^{2}}
\end{equation*}
We measured the errors for various values of $\delta$, taking $\sigma=1$. The results are plotted in Figure \ref{fig:schroErrors}. Note further reducing $\delta$ does not reduce the error significantly beyond $10^{-8}$; we believe this is due to floating point errors.

\subsection{Euler equations, linearized around a jet}

\label{sec:Euler}

We now move on to a more interesting example.

We consider now the two-dimensional Euler equations linearized around a homogeneous jet flow in the $x_{1}$ direction. The vector $\vu(x,t)$ is 3 dimensional, with $\vu_{1}$ representing the pressure change, $\vu_{2,3}$ representing the velocity field in the $x_{1,2}$ directions (respectively). In this case, $\HH$ takes the form:
\begin{equation}
  \label{eq:linearizedEuler}
  \HH =  \left[
    \begin{array}{lll}
      M \partial_{x_{1}} & - \partial_{x_{1}} & -\partial_{x_{2}}\\
      - \partial_{x_{1}} & M \partial_{x_{1}} & 0\\
      - \partial_{x_{2}} & 0 & M \partial_{x_{1}} \\
    \end{array}
  \right]
\end{equation}
Here, $0 \leq M < 1$ is the mach number. $\HH$ has eigenvalues $\omega_{1}(\vk) = M k_{1}+\absSmall{\vk}$, $\omega_{2}(\vk) = M k_{1}-\absSmall{\vk}$ and $\omega_{3}(\vk) = M k_{1}$. The diagonalizing matrix is
\begin{equation}
  D = \frac{1}{\sqrt{2}\absSmall{\vk}}
  \left[
    \begin{array}{lll}
      - \absSmall{\vk} & k_{1}& k_{2}\\
      \absSmall{\vk} & k_{1} & k_{2}\\
      0  &  -\sqrt{2}k_{1} & \sqrt{2}k_{2}
    \end{array}
  \right].
\end{equation}

When $M \neq 0$, this is the classical example of an equation for which the PML becomes unstable \cite{hu:unstablePML}.

\begin{figure}
\setlength{\unitlength}{0.240900pt}
\ifx\plotpoint\undefined\newsavebox{\plotpoint}\fi
\sbox{\plotpoint}{\rule[-0.200pt]{0.400pt}{0.400pt}}%
\includegraphics[scale=0.6]{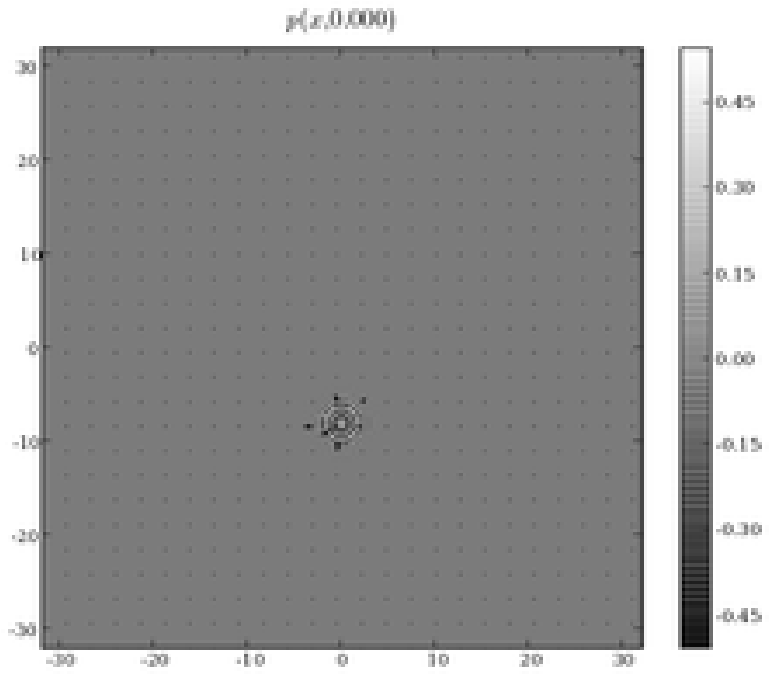}
\includegraphics[scale=0.6]{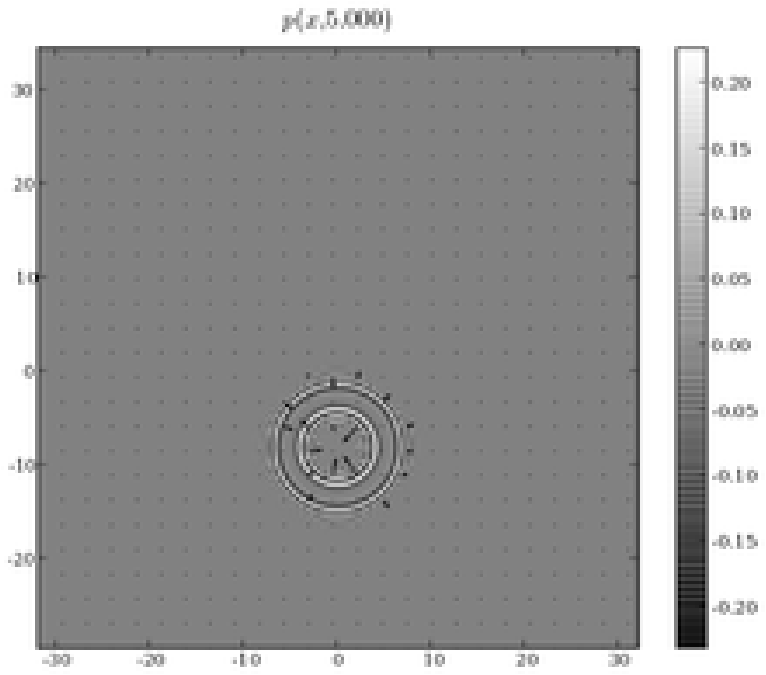}
\includegraphics[scale=0.6]{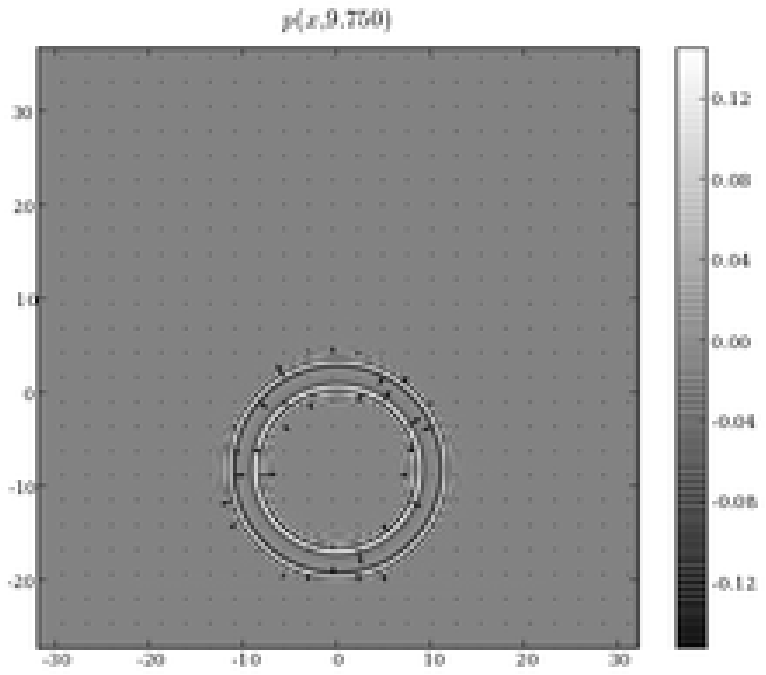}
\includegraphics[scale=0.6]{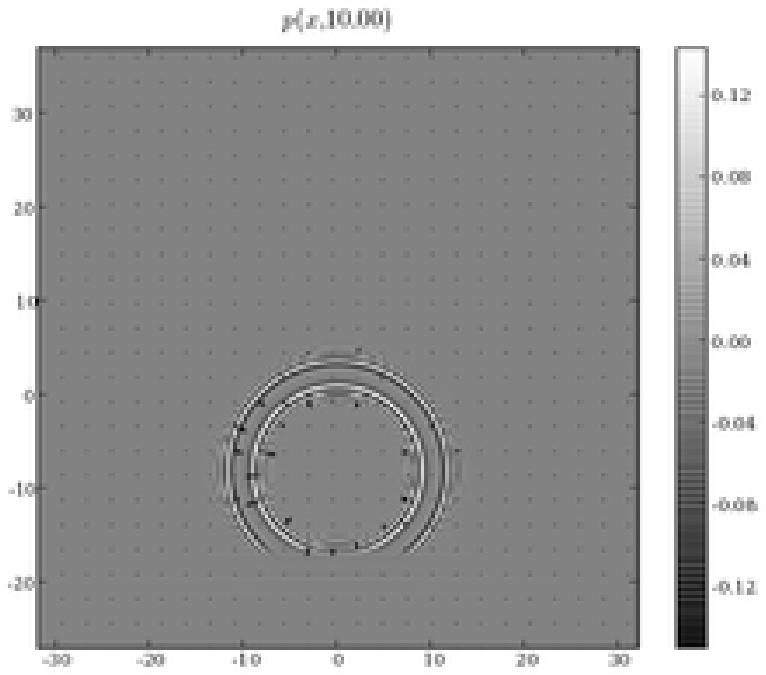}
\includegraphics[scale=0.6]{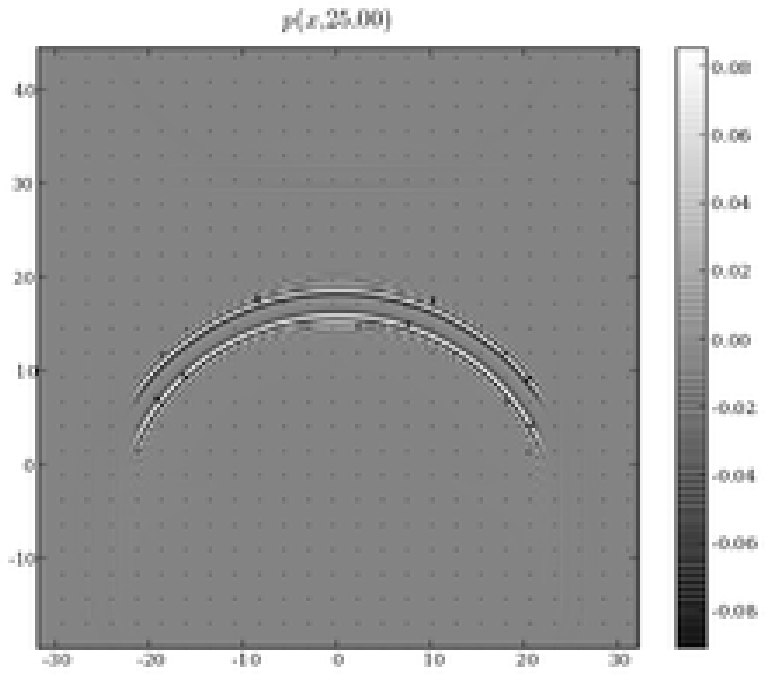}
\includegraphics[scale=0.6]{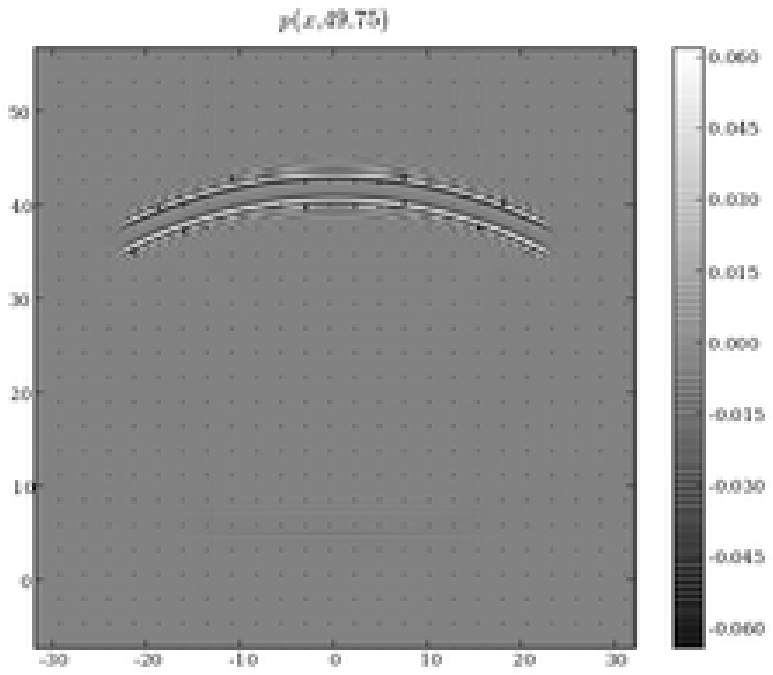}

\caption{The pressure field at various times. The arrows indicate the direction of the velocity field. Note that the phase space filter was applied between $t=9.75$ and $t=10.0$. This simulation is available in movie form from http://cims.nyu.edu/$\sim$stucchio/software/kitty/demos/euler\_tdpsf.mpg. }
\label{fig:eulerFrames1}
\end{figure}

\begin{figure}
\setlength{\unitlength}{0.240900pt}
\ifx\plotpoint\undefined\newsavebox{\plotpoint}\fi
\sbox{\plotpoint}{\rule[-0.200pt]{0.400pt}{0.400pt}}%
\includegraphics[scale=0.5]{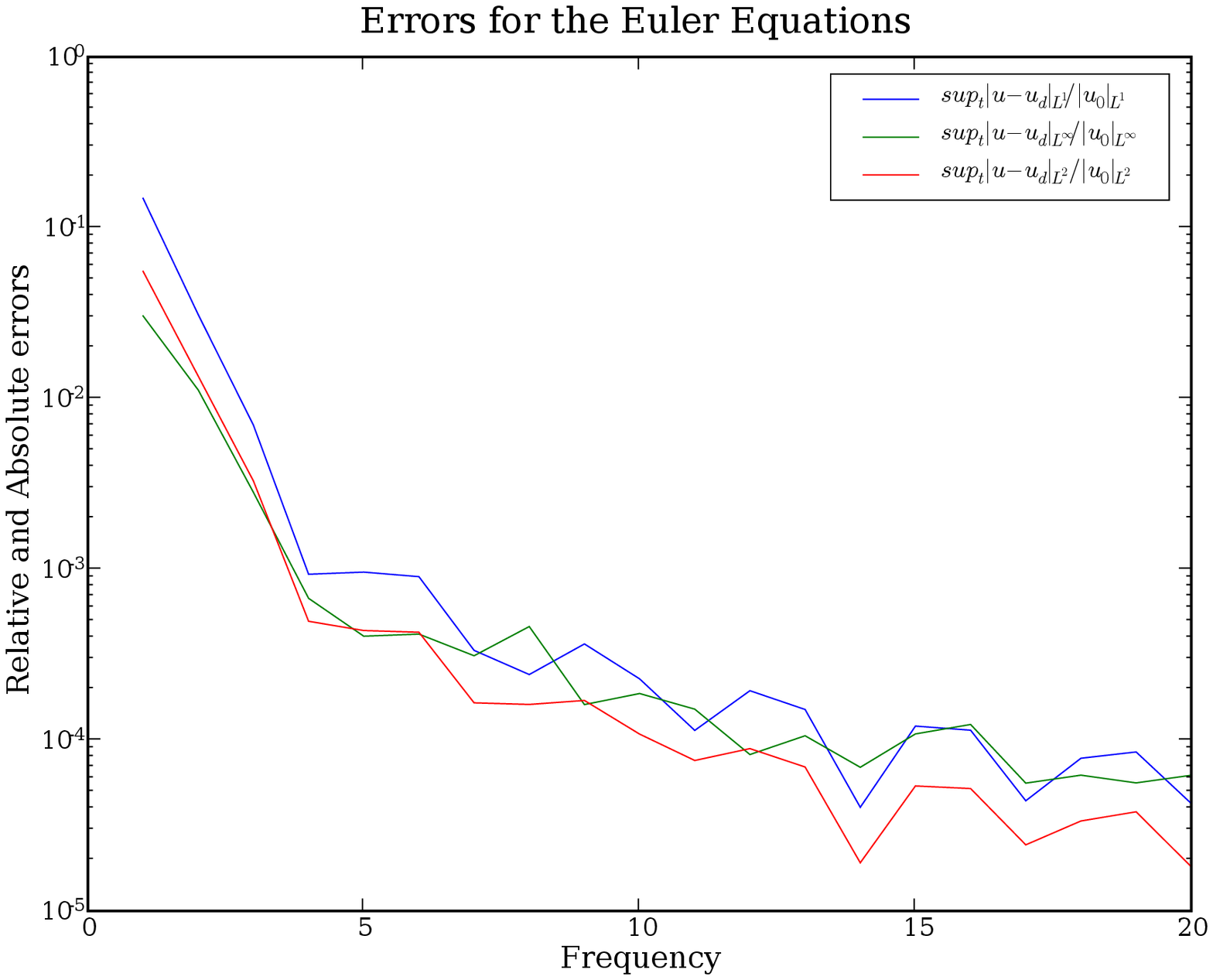}

\caption{The relative error (measured in various norms) as a function of the frequency of the initial condition.}
\label{fig:eulerErrors}
\end{figure}

We solved the Euler equations on the region $[-32,32]^{2}$ with lattice spacing $\delta x = 0.125$ (for a total of $512^{2}$ lattice points). Such a setup is valid for spatial frequencies up to $\kmax = \pi/\delta x = 25.1$. We solved the system using the Fast Fourier transform to compute \eqref{eq:propagatorDef}, which is accurate to machine precision provided no waves reach the boundary (\cite[Theorem 4.1]{us:TDPSFrigorous} or \cite[Theorem A.1]{us:multiscale}). This accuracy is independent of time-step, which we took to be $\delta t = 0.25$ in order to get a watchable frame-rate in the plots. The phase space filter region was taken to have width $16$ ($128$ lattice points), with $\Tstep=1.5$, and $\sigma=1.0$. The initial condition was taken to be

\begin{subequations}
  \begin{equation}
    \label{eq:eulerInitialCondition}
    u_{1}(x,t=0) = r^{2} e^{-r^{2}/9} \cos(K r)
  \end{equation}
  with
  \begin{equation}
    r = \sqrt{(x-8)^{2}+y^{2} }.
  \end{equation}
\end{subequations}
with the $K$ varying from $1$ to $20$ (the frequency range of the problem). This yields an initial condition with frequency localized near $\absSmall{\vk}=K$. The results were compared to another simulation on a large box. The result is that for $K > 4$, $L^{2}$ error of $10^{-3}$ (relative to the initial condition) is achieved up to $t=50$ (see Figure \ref{fig:eulerErrors}). The error at low frequencies can be dealt with by various means, see Section \ref{sec:lowFrequency} and \cite{us:multiscale}.

\subsection{Maxwell's Equations in an Orthotropic Medium}

\label{sec:Maxwell}

Let $\vec{E}$ be the electric field and $\vec{H}=\mu \vec{B}$ with $\vec{B}$ the magnetic field. Let $\epsilon$ and $\mu$ be the electrical and magnetic permeability's of a medium.

In order to bring Maxwell's equations to symmetric form, we introduce the auxiliary variable $\vu = ( \sqrt{\mu} \vec{H}, \sqrt{\epsilon} \vec{E})^{T}$. With this formulation, the Hamiltonian for Maxwell's equations can be written in block form:

\begin{equation}
  \HH =
  \left[
    \begin{array}{ll}
      0 &  -\mu^{-1/2} \nabla \times \epsilon^{-1/2} \\
      \epsilon^{-1/2} \nabla \times \mu^{-1/2}  & 0\\
    \end{array}
  \right]
\end{equation}
where $\nabla \times$ is interpreted as a matrix,
\begin{equation*}
  \nabla \times =
  \left[
    \begin{array}{lll}
      0 & -\partial_{z} & \partial_{y} \\
      \partial_{z} & 0 & -\partial_{x} \\
      - \partial_{y} & \partial_{x} & 0
    \end{array}
  \right].
\end{equation*}
The symmetry of $\epsilon$, $\mu$ and $\nabla \times$ implies that $\HH$ is skew-adjoint.

To further simplify the system, we make the following additional assumptions. We assume $\mu=1$. We assume the system is $z$-independent, i.e. $\partial_{z} = 0$, and we can restrict ourselves to two-dimensional simulations. Lastly, we simplify the electrical permittivity:
\begin{equation}
  \epsilon =
  \left[
    \begin{array}{lll}
      1 & b & 0\\
      b & 1 & 0\\
      0 & 0 & c
    \end{array}
  \right]
\end{equation}
This is the simplest possible birefringent system. With the variables
\begin{eqnarray*}
  f &=& (1/2)(\sqrt{1+b}+\sqrt{1-b}) \\
  g &=& (1/2)(-\sqrt{1+b}+\sqrt{1-b}),
\end{eqnarray*}
we can write the dispersion relation as:
\begin{subequations}
  \begin{eqnarray}
    \omega_{j=1,2}(\vk) & =& (-1)^{1+j} i c^{-1} \absSmall{\vk}\\
    \omega_{j=3,4}(\vk) & = & (-1)^{1+j} i \sqrt{(f^{2}+g^{2})(k_{1}^{2}+k_{2}^{2}) -4fg k_{1}k_{2} }\\
    \omega_{j=5,6}(\vk) & = & 0.
  \end{eqnarray}
\end{subequations}
with
\begin{equation}
  D = \left[
    \begin{array}{llllll}
      \frac{- k_{2}}{\sqrt{2}\abs{\vk}} &\frac{k_{1}}{\sqrt{2}\abs{\vk}} & 0 & 0 & 0 & 2^{-1/2} \\
      \frac{k_{2}}{\sqrt{2}\abs{\vk}}& \frac{-k_{1}}{\sqrt{2}\abs{\vk}}& 0 & 0 & 0 & 2^{-1/2} \\
      0 & 0 & -2^{-1/2} & -\frac{f k_{2} - g k_{1}}{\sqrt{2} E(\vk)} & \frac{f k_{1} - g k_{2}}{\sqrt{2} E(\vk)} & 0 \\
      0 & 0 & 2^{-1/2} & -\frac{f k_{2} - g k_{1}}{\sqrt{2} E(\vk)} & \frac{f k_{1} - g k_{2}}{\sqrt{2} E(\vk)} & 0 \\
     k_{1} \abs{\vk}^{-1} & k_{2} \abs{\vk}^{-1} & 0 & 0 & 0 & 0 \\
     0 & 0 & 0 & \frac{f k_{1} - g k_{2}}{E(\vk)} & \frac{f k_{2} - g k_{1}}{E(\vk)} & 0
    \end{array}
  \right]
\end{equation}
with $E(\vk)=\sqrt{(f^{2}+g^{2})(k_{1}^{2}+k_{2}^{2}) -4fg k_{1}k_{2} }$.

This can be further simplified by noting that $u_{1},u_{2}$ and $u_{6}$ (corresponding to $B_{x},B_{y},E_{z}$ are uncoupled to $u_{3},u_{4}$ and $u_{5}$. The $u_{1}, u_{2}$ and $u_{6}$ modes (corresponding to $B_{x},B_{y}$ and $E_{z}$) are actually isotropic waves, so we will ignore them. We therefore restrict consideration to $u_{3}, u_{4}$ and $u_{5}$ ($u_{3}$ corresponds to $B_{z}$, and $u_{4,5}$ correspond to to linear combinations of $E_{x}$ and $E_{y}$).

\begin{figure}
\setlength{\unitlength}{0.240900pt}
\ifx\plotpoint\undefined\newsavebox{\plotpoint}\fi
\sbox{\plotpoint}{\rule[-0.200pt]{0.400pt}{0.400pt}}%
\includegraphics[scale=0.6]{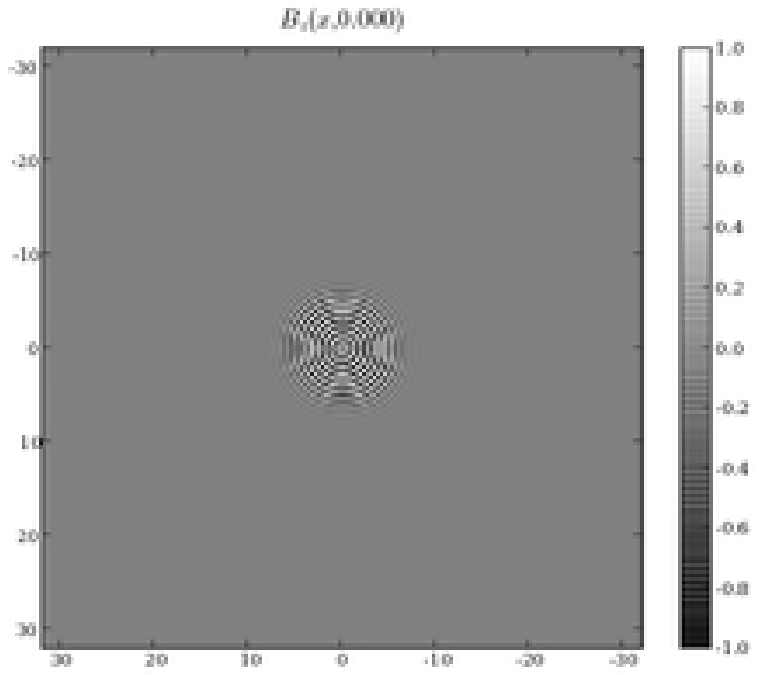}
\includegraphics[scale=0.6]{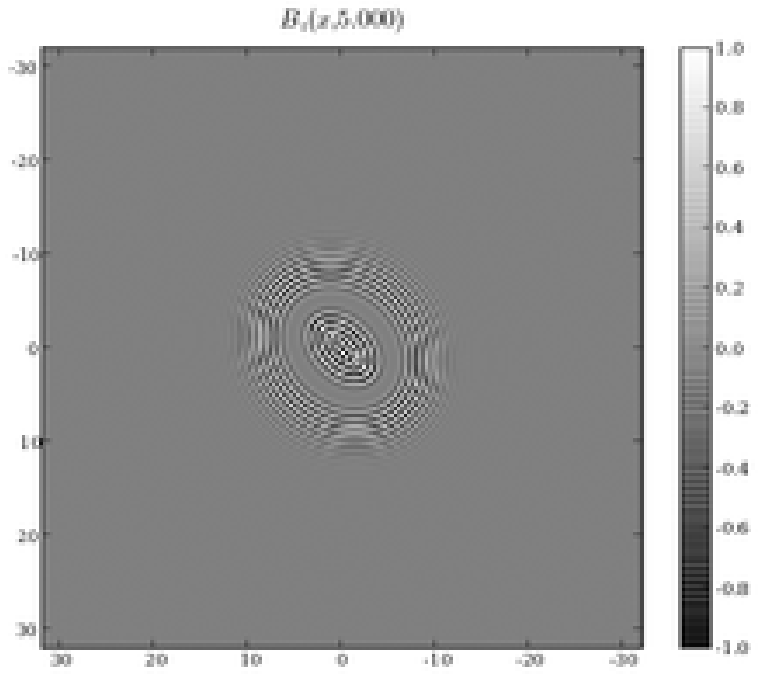}
\includegraphics[scale=0.6]{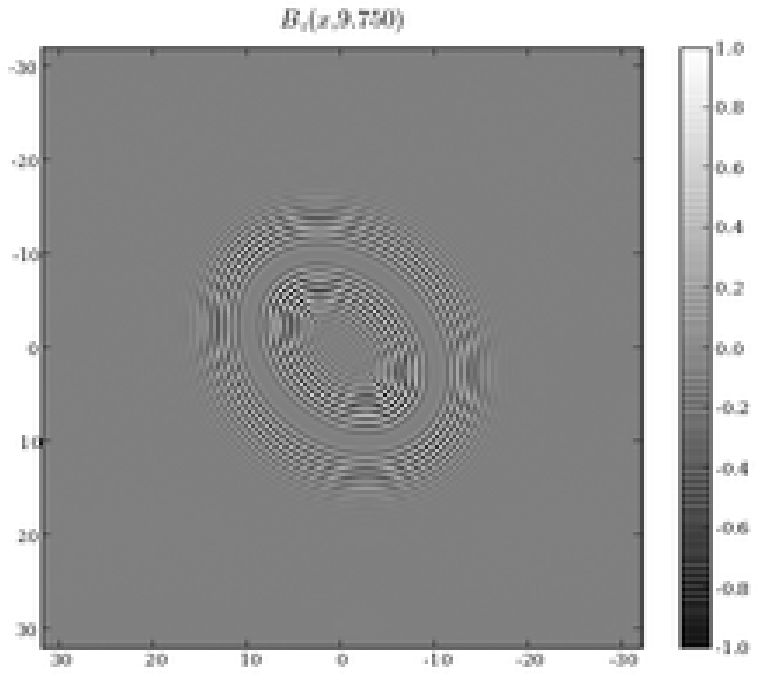}
\includegraphics[scale=0.6]{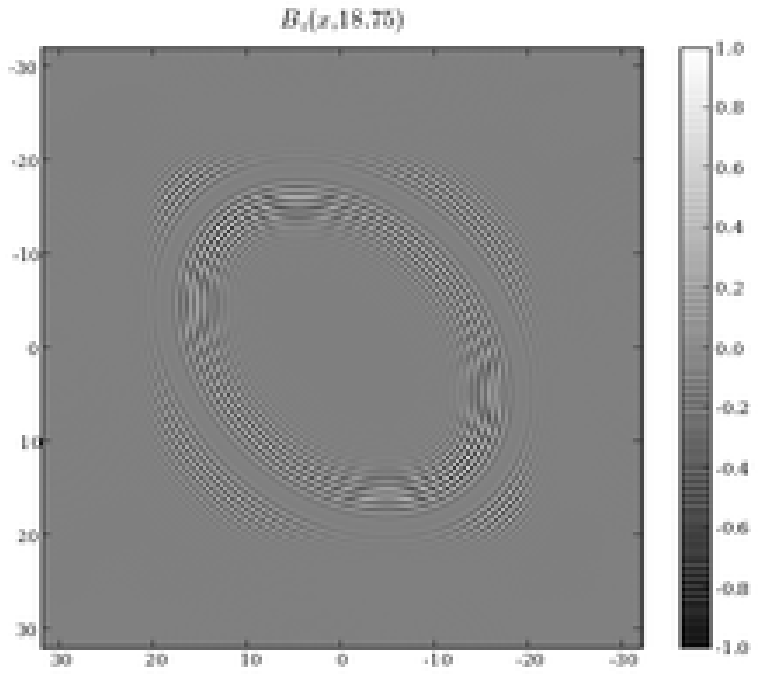}
\includegraphics[scale=0.6]{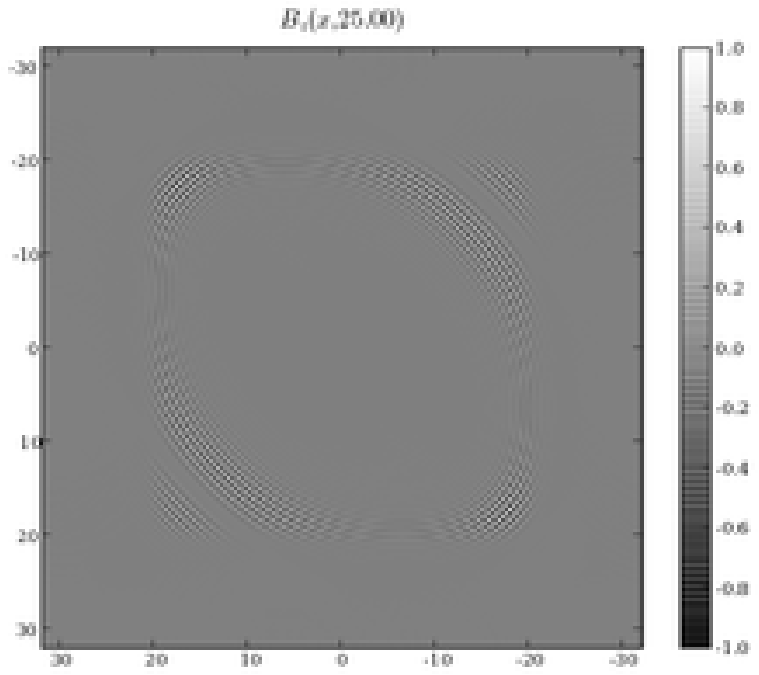}
\includegraphics[scale=0.6]{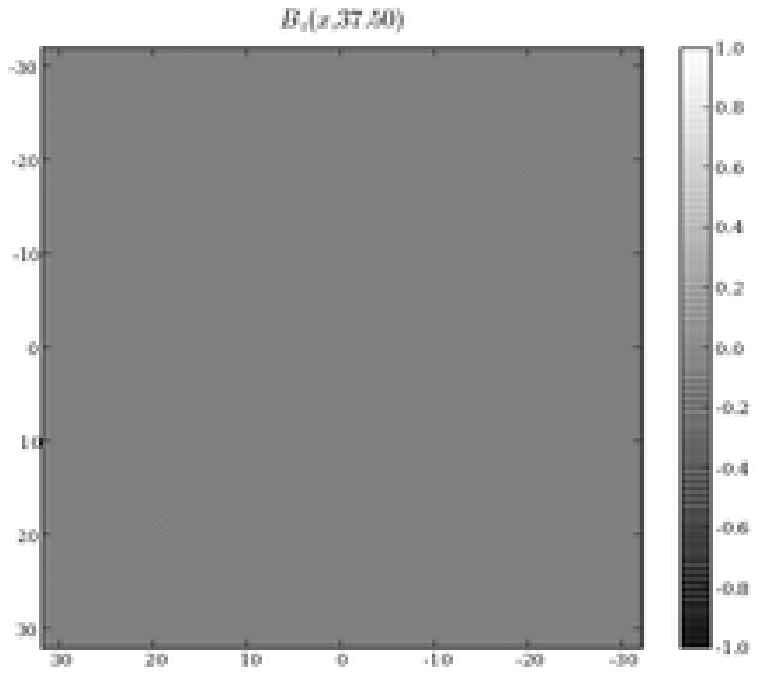}

\caption{The magnetic component $B_{z}$ of the electromagnetic pulse at various times. Note that the phase space filter was applied between $t=9.75$ and $t=10.0$. The non-radial shape of the wave is due to the anisotropy of the medium. This simulation is available in movie form from http://cims.nyu.edu/$\sim$stucchio/software/kitty/demos/maxwell\_tdpsf.mpg. }
\label{fig:maxwellFrames1}
\end{figure}

\begin{figure}
\setlength{\unitlength}{0.240900pt}
\ifx\plotpoint\undefined\newsavebox{\plotpoint}\fi
\sbox{\plotpoint}{\rule[-0.200pt]{0.400pt}{0.400pt}}%
\includegraphics[scale=0.5]{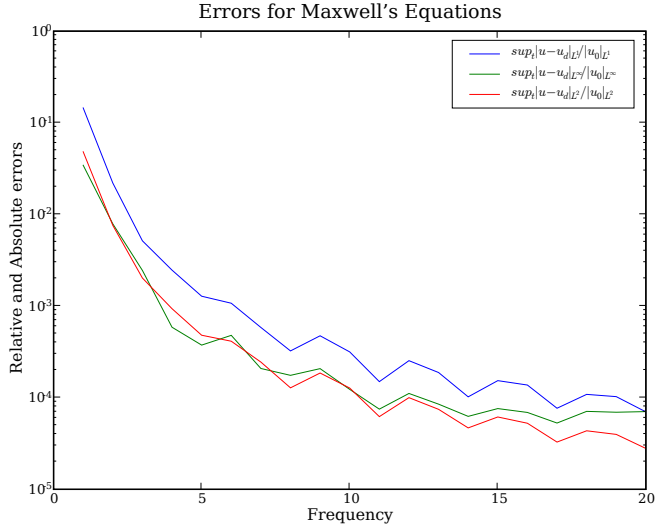}

\caption{The relative error (measured in various norms) as a function of the frequency of the initial condition.}
\label{fig:maxwellErrors}
\end{figure}

We solved the reduced system of Maxwell's equations with the same parameters as in Section \ref{sec:Euler}, with the same initial condition (replacing $u_{1}$ by $u_{3}$ in \eqref{eq:eulerInitialCondition}). The anisotropy parameter $b$ was chosen to be $0.25$. The results are comparable to those for the Euler equation, see Figures \ref{fig:maxwellFrames1} and \ref{fig:maxwellErrors}.

\subsection{Long Time Stability}

In Theorem \ref{thm:stability} of Section \ref{sec:stability}, it was proved that the phase space filtering algorithm is stable. To demonstrate the validity of the theorem, we ran a simulation of the Euler equations and Maxwell's equations up to time $t=2000$. While we can not determine the accuracy over such long time intervals (the reference simulation would require an extremely large box), we can study the growth of the $L^{2}$-norm.

\begin{figure}
\setlength{\unitlength}{0.240900pt}
\ifx\plotpoint\undefined\newsavebox{\plotpoint}\fi
\sbox{\plotpoint}{\rule[-0.200pt]{0.400pt}{0.400pt}}%
\includegraphics[scale=0.5]{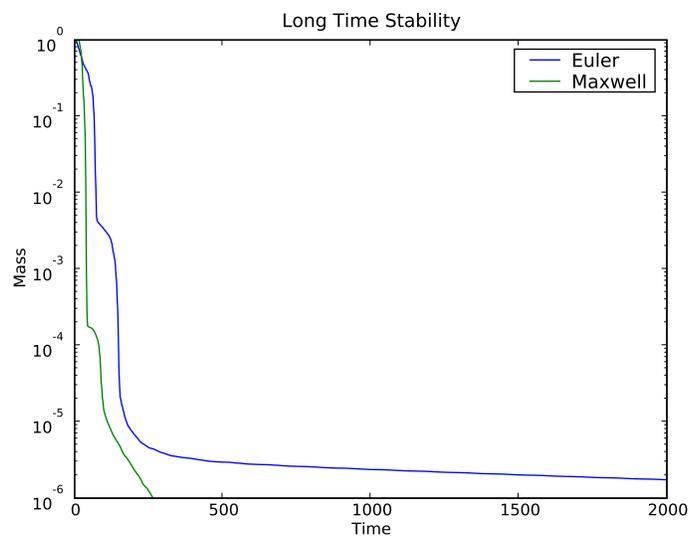}

\caption{The $L^{2}$ norm of solutions of the Euler and Maxwell equations as a function of time. In both cases, we took $K=10$ (recall \eqref{eq:eulerInitialCondition}). The Mach number was $M=0.5$ for the Euler equations and the anisotropy was $b=0.25$ for Maxwell's equations.}
\label{fig:longTimeStability}
\end{figure}

The results of such a simulation are plotted in Figure \ref{fig:longTimeStability}. They indicate that Theorem \ref{thm:stability} is correct, and that the mass of the solution decreases monotonically with time.

\section{The Low Frequency Problem}
\label{sec:lowFrequency}

As is apparent from Figures \ref{fig:eulerErrors} and \ref{fig:maxwellErrors}, the phase space filtering approach does not work well for waves with low frequency. The reason for this is that to localize in frequency, the region in which one works must be $O(1)$ wavelengths long. The simplest remedy is to increase the width of the filter. If the smallest frequency relevant to the problem is $\kb$, then the width of the buffer is $O(\kb^{-1})$, which means that the computational cost is of order $O(\kb^{-\ndom})$.

This problem can be remedied by a somewhat more involved method, which has been implemented for the Schr\"odinger equation \cite{us:multiscale}. The essential idea is to increase the width of the box, but reduce the sampling rate on the extended region. Then high frequency waves are filtered at the edge of the highly sampled region, and low frequency waves are filtered on the edge of the coarsely sampled region, but using a wider filter capable of resolving low frequency waves.

With this method, even though the box has width $O(\kb^{-1})$, the number of samples required is only $O(\log(\kmax / \kb))$ (computation time scales similarly, up to logarithmic prefactors). This allows resolution of outgoing waves at low frequencies in logarithmic rather than linear cost. It is also argued heuristically in \cite{us:multiscale} that this computational complexity is close to the best possible.

{\bf Acknowledgements:} We thank P. Petropoulos for pointing out to us the problem of PML instability, and R. Goodman for a useful discussion on lattice waves.

\nocite{MR966733}
\nocite{MR0436612}
\nocite{MR0471386}
\nocite{MR517938}

\nocite{MR1819643}
\nocite{MR2032866}

\bibliographystyle{hplain}
\bibliography{../stucchio.bib}

\end{document}